\documentclass[dvips,12pt]{amsart2000}
\usepackage[dvips]{graphicx}
\usepackage[usenames,dvipsnames]{color}
\usepackage{amsmath2000,amsthm2000,amssymb,a4wide,times}
\usepackage[mathscr]{eucal}

\newcommand{\dddots}{\mathinner{\mkern1mu\raise1pt\vbox{\kern7pt\hbox{.}}\mkern2mu\raise4pt\hbox{.}\mkern2mu\raise7pt\hbox{.}\mkern1mu}}

\newtheorem{theorem}{Theorem}[section]
\newtheorem{lemma}[theorem]{Lemma}

\advance\textheight by 8pt

\begin{document}

\title{Boundary Representations for Operator Algebras}
\author[M.~A.~Dritschel and S.~A.~McCullough]{Michael A. Dritschel and Scott McCullough${}^*$}

\address{Department of Mathematics, 
School of Mathematics and Statistics,
Merz Court, University of Newcastle upon Tyne,
Newcastle upon Tyne,
NE1 7RU
UK}

\email{M.A.Dritschel@newcastle.ac.uk}

\address{Department of Mathematics,
University of Florida,
Box 118105,
Gainesville, FL 32611-8105
USA}

\email{sam@math.ufl.edu}



\thanks{${}^*$research supported by the NSF}

\begin{abstract}
  All operator algebras have (not necessarily irreducible) boundary
  representations.  A unital operator algebra has enough such boundary
  representations to generate its $C^*$-envelope.
\end{abstract}

\maketitle

\section{Introduction}
\label{sec:introduction}

Concretely, an operator algebra $\mathcal{A}$ is a subalgebra of
$\mathcal{B}(K)$, the bounded linear operators on some Hilbert space
$K$.  It is unital if it contains the identity operator.  The algebra
$M_\ell(\mathcal{A})$ of $\ell \times \ell$ matrices with entries from
$\mathcal{A}$ inherits a norm as a subspace of
$M_\ell(\mathcal{B}(K))$ identified canonically with
$\mathcal{B}(\oplus_1^\ell K)$.  The Blecher, Ruan and Sinclair
Theorem \cite{MR91b:47098} characterizes unital operator algebras in
terms of a matrix norm structure, while a theorem of Blecher
\cite{MR96k:46098} does the same for non-unital algebras assuming the
algebra multiplication is completely bounded.  Consequently it is
possible to speak abstractly of an operator algebra without reference
to an ambient $\mathcal{B}(K)$.

A linear mapping $\phi:\mathcal{A}\to\mathcal{B}(H)$ induces a
linear mapping $\phi_\ell:M_\ell(\mathcal{A}) \to
\mathcal{B}(\oplus_1^\ell H)$ by applying $\phi$ entry-wise, so that
$\phi_\ell (a_{jm}) =(\phi(a_{jm}))$.  The map $\phi$ is completely
bounded if $\phi$ is bounded and there exists $C$, independent of
$\ell$, such that $\|\phi_\ell\|\leq C$, it is completely contractive
if it is completely bounded with $C \leq 1$, and it is completely
isometric if $\phi_\ell$ is an isometry for each $\ell$.  Finally, a
representation of $\mathcal{A}$ on the Hilbert space $H$ is an algebra
homomorphism $\phi:\mathcal{A}\to\mathcal{B}(H)$.  If
$\mathcal{A}$ is unital, it is assumed that any representation of
$\mathcal{A}$ takes the unit to the identity operator.

A \textit{boundary representation} (\cite{MR52:15035},
\cite{MR40:6274}) of the unital operator algebra $\mathcal{A}$
consists of a homomorphism $\phi:\mathcal{A}\to\mathcal{C}$, where
$\mathcal{C}$ is a $C^*$-algebra and $C^*(\phi(\mathcal{A}))
=\mathcal{C},$ together with a representation $\pi:\mathcal{C}\to
\mathcal{B}(H)$ such that the only completely positive map on
$\mathcal{C}$ agreeing with $\pi$ on $\phi(\mathcal{A})$ is $\pi$
itself.  In originally defining boundary representations, Arveson also
required that they be irreducible.  We do not impose this condition.

The $C^*$-\textit{envelope} of $\mathcal{A}$, denoted
$C^*_e(\mathcal{A})$, is the essentially unique smallest $C^*$-algebra
amongst those $C^*$-algebras $\mathcal{C}$ for which there is a
completely isometric homomorphism $\phi:\mathcal{A}\to \mathcal{C}$.
For instance, if $\mathcal{A}$ is a uniform algebra, then
$C^*_e(\mathcal{A})$ is the $C^*$-algebra of continuous functions on
the {\v S}ilov boundary of $\mathcal{A}$.  In fact, in this case the
irreducible boundary representations correspond to peak points of
$\mathcal{A}$.  Arveson proved that $C^*_e(\mathcal{A})$ exists
provided there are enough boundary representations for $\mathcal{A}$.
However, the existence of $C^*_e(\mathcal{A})$ does not imply the
existence of boundary representations and Hamana \cite{MR81h:46071}
established the existence of $C^*_e(\mathcal{A})$ in general without
recourse to boundary representations.

In this note we show, by elaborating on a construction of Agler
essential to his approach to model theory \cite{MR90a:47016} and using
a characterization of boundary representations due to Muhly and Solel
\cite{MR99f:47057}, that boundary representations exist, and then
following an argument similar to Arveson's, we also derive the
existence of $C^*_e(\mathcal{A})$.

Agler's approach is to consider a \textit{family}
$\mathcal{F}_{\mathcal{A}}$ of representations of an algebra
$\mathcal{A}$ (which is not necessarily an operator algebra).  This is
a collection of representations which is
\begin{enumerate}
\item Closed with respect to direct sums (so if $\{\pi_\alpha\}$ is an
  arbitrary set of representations in the family, then
  $\bigoplus_\alpha \pi_\alpha$ is also a representation in the
  family);
\item Hereditary (that is, if $\phi$ is a representation in the family
  and $L$ is a subspace which is invariant for all $\phi(a)$, $a\in
  \mathcal{A}$, then $\phi|_L$, the restriction of $\phi$ to $L$, is
  also in the family);
\item Closed with respect to unital $*$-representations (so if $\phi:
  \mathcal{A} \to \mathcal{B}(H)$, and $\nu:\mathcal{B}(H) \to
  \mathcal{B}(K)$ is a unital $*$-representation, then $\nu\circ\phi$
  is in the family).

\smallskip

\noindent\hskip-\leftmargin If $\mathcal{A}$ is non-unital, then we
also require

\smallskip

\item $\mathcal{A}$ is closed with respect to spanning representations
  with respect to the partial ordering on dilations (defined below).
  A consequence of (1) is that the norms of all $\pi\in
  \mathcal{F}_{\mathcal{A}}$ are uniformly bounded.
\end{enumerate}
When $\mathcal{A}$ is unital, (1)--(3) can be shown to imply (4),
using an argument similar to that used to prove (1) of Theorem 1.1 in
\cite{MR2000b:47025}.

Special examples of families include the collection of all completely
contractive representations of an operator algebra $\mathcal{A}$ and
the collection of all representations $\pi$ of the disc algebra such
that $\pi(z)$ is an isometry.

A representation $\phi$ \textit{lifts} to a representation $\psi$ if
$\phi$ is the restriction of $\psi$ to an invariant subspace.
Following Agler \cite{MR90a:47016} we will say that the representation
$\phi:\mathcal{A}\to \mathcal{B}(H)$ is \textit{extremal}, if whenever
$K$ is a Hilbert space containing $H$ and $\psi:\mathcal{A} \to B(K)$
is a representation such that $H$ is invariant for $\psi(\mathcal{A})$
and $\phi=\psi|_H$, then $H$ reduces $\psi(\mathcal{A})$.  Further, if
$\rho:\mathcal{A}\to \mathcal{B}(L)$ is a representation, then $\rho$
lifts to an extremal representation; i.e., there exists a Hilbert
space $H$ containing $L$ and an extremal representation $\phi$ of
$\mathcal{A}$ such that $L$ is invariant for $\phi(\mathcal{A})$ and
$\rho=\phi|_L$ (\cite{MR90a:47016}, Proposition 5.10).

Lifting induces a partial ordering on representations, with
$\phi_\alpha\leq {\phi_\beta}$ being equivalent to ${\phi_\alpha}$
lifting to ${\phi_\beta}$.  If $\mathcal{S}$ is a totally ordered set
of liftings (with respect to this partial ordering), then we define
the \textit{spanning representation} $\phi_s:\mathcal{A}\to
\mathcal{B}(H_s)$ by setting ${H}_s$ to be the closed span of the
${H}_\alpha$'s over all ${\alpha} \in \mathcal{S}$, and then densely
defining $\phi_s$ to be $\phi_\alpha$ on $H_\alpha$ and extending to
all of $H_s$ by boundedness of the representations $\phi_{\alpha}$.
It is readily verified that $\phi_s$ is a representation which lifts
each $\phi_\alpha$.

The representation $\phi:\mathcal{A}\to \mathcal{B}(H)$ dilates to the
representation $\psi:\mathcal{A}\to \mathcal{B}(K)$ if $K$ contains
$H$ and $\phi(a)=P_H\psi(a)|_H$ for all $a\in\mathcal{A}$.  A
fundamental result of Sarason \cite{MR32:6229} says that a
representation $\phi$ dilates to a representation $\psi$ if and only
if $H$ is semi-invariant for $\psi$.  Thus, there exists subspaces
$L\subset N\subseteq K$ invariant for $\psi$ such that $H=N\ominus L$.
Alternatively, $K = L\oplus H \oplus M$ with $L$ and $L\oplus H$
invariant for $\phi$.  Just as in the case of liftings, dilating
induces a partial ordering on representations in the obvious manner.
We can also similarly define spanning representations of totally
ordered sets of representations, and this is what is used in item (4)
above.  Note that liftings are also dilations (with $L=\{0\}$).  Hence
the partial ordering on dilations subsumes that of liftings, and in
particular, any spanning representation of liftings is one in terms of
dilations as well.

As was noted above, families of representations over unital algebras
contain all spanning representations formed from chains of
representations in the family, though it appears that this is a needed
added assumption in the non-unital setting.  On the other hand, there
are interesting collections of representations which are closed with
respect to (1) and (4) of a family, but not necessarily (2) and (3).
For example, the collection of all completely isometric
representations of an algebra fall into this category.  Since the
theorems we prove below only depend on existence of spanning
representations in the collections of representations we are
considering and all representations being uniformly bounded, we define
the \textit{extended family} of an algebra $\mathcal{A}$ to be a
collection of representations of $\mathcal{A}$ which is closed under
the formation of direct sums and spanning representations.

For dilations, the equivalent of an extremal will be referred to as a
$\partial$-\textit{representation}.  The representation
$\phi:\mathcal{A}\to \mathcal{B}(H)$ is a $\partial$-representation if
whenever $\psi:\mathcal{A} \to B(K)$ dilates $\phi$, then $H$ reduces
$\psi(\mathcal{A})$.

Muhly and Solel \cite{MR99f:47057} show, in the language of Hilbert
modules rather than representations, that for unital operator algebras
$\partial$-representations coincide with boundary representations
(forgetting the irreducibility requirement).

\begin{theorem}\label{th:0}
  Let $\mathcal{A}$ be a unital operator algebra.  Then
  $\rho:\mathcal{A} \to \mathcal{B}(H)$ is a $\partial$-representation
  if, and only if, given any completely isometric map $\phi :
  \mathcal{A} \to \mathcal{C}$ where $\mathcal{C}$ is a $C^*$ algebra
  with $\mathcal{C} = C^*(\phi(\mathcal{A}))$, there exists a
  boundary representation $\pi:C^*(\rho(\mathcal{A})) \to
  \mathcal{B}(H)$ such that $\pi \circ \phi = \rho$.
\end{theorem}

The proof of Muhly and Solel of this equivalence uses the existence of
the $C^*$-envelope.  Our main result and proof of the existence of the
$C^*$-envelope do not depend on their work.  However, it should be
noted that a proof of the equivalence which does not already assume
the existence of the $C^*$-envelope is possible, and we sketch a proof
below based along a line of reasoning in (\cite{McSPV}, Theorem 1.2).

\begin{proof}[Sketch of the proof of (\ref{th:0})]
  Suppose $\phi : \mathcal{A} \to \mathcal{C} =
  C^*(\phi(\mathcal{A}))$ is completely isometric and
  $\pi:C^*(\phi(\mathcal{A})) \to \mathcal{B}(H)$ is a boundary
  representation.  Set $\rho = \pi \circ \phi$, and note that it is
  completely contractive.  Suppose $\nu:\mathcal{A} \to
  \mathcal{B}(K)$ dilates $\rho$.  The goal is to show that $H$
  reduces $\nu$.
  
  To this end, define a map $\gamma:\phi(\mathcal{A}) \to
  \mathcal{B}(K)$ by $\gamma(\phi(a)) = \nu(a)$, $a\in\mathcal{A}$.
  This map is completely contractive, and so by the Arveson extension
  theorem extends to a completely positive unital map
  $\gamma:C^*(\phi(\mathcal{A})) \to \mathcal{B}(K)$ with
  $\gamma\circ\phi = \nu$.  Observe that the map which takes $b
  \mapsto P_{H} \gamma(\phi(b)) |_H$, $b\in C^*(\phi(\mathcal{A}))$ is
  completely positive, and by definition, $P_{H}
  \gamma(\phi(a)) |_H = \rho(a) = \pi(\phi(a))$ for all $a\in
  \mathcal{A}$.  We have assumed that $\pi$ is a boundary
  representation, so in fact $P_{H} \gamma(\phi(b)) |_H = \pi(b)$ for
  all $b\in C^*(\phi(\mathcal{A}))$.  From this we have for all $a\in
  \mathcal{A}$,
  \begin{eqnarray*}
    \rho(a)\rho(a)^* &=& \pi(\phi(a))\pi(\phi(a))^* \\
    &=& \pi(\phi(a)\phi(a)^*) \\
    &=& P_H \gamma(\phi(a)\phi(a)^*) |_H \\
    &\geq& P_H \gamma(\phi(a))\gamma(\phi(a)^*) |_H = P_H
    \nu(a)\nu(a)^* |_H \\
    &\geq& P_H \gamma(\phi(a)) P_H \gamma(\phi(a)^*) |_H
    = P_H \nu(a) P_H \nu(a)^* |_H \\
    &=& \rho(a)\rho(a)^*,
  \end{eqnarray*}
  where the first inequality is the Cauchy Schwarz inequality for
  completely positive maps \cite{MR88h:46111}.  From this we see that
  $\nu(a)^* H \subseteq H$.  An identical argument gives $\nu(a) H
  \subseteq H$, proving that $H$ reduces $\nu$.

  The converse is a straightforward exercise and is left to the
  reader.  
\end{proof}

In the course of this paper, we will prove the following.

\begin{theorem}\label{th:1}
  If $\rho:\mathcal{A}\to \mathcal{B}(H)$ is a representation in an
  extended family $\mathcal{F}_{\mathcal{A}}$, then there exists a
  Hilbert space $K$ containing $H$ and a $\partial$-representation
  $\phi:\mathcal{A}\to \mathcal{B}(K)$ also in
  $\mathcal{F}_{\mathcal{A}}$ such that $\rho$ dilates to $\phi$.
\end{theorem}

As mentioned above, Arveson's original definition of boundary
representation required $\pi$ to be irreducible.  Note that
(\ref{th:1}) does not imply the existence of irreducible boundary
representations.

The authors would like to thank Jim Agler for sharing a draft of the
first several chapters of a recent manuscript on his abstract model
theory and Vern Paulsen for his valuable assistance and an advance
look at the new edition of his book, \textsl{Completely Bounded Maps
  and Dilations}, particularly the chapters on the abstract
characterization of operator algebras and the $C^*$-envelope.

The remainder of the paper is organized as follows.  Section
\ref{sec:lift-extr-agler} establishes (\ref{th:2}) giving the
existence of extremals, which form the core of any model in Agler's
approach to model theory \cite{MR90a:47016}.  Although versions of
this result are quite old, no proofs yet appear in the literature.  In
Section \ref{sec:dilat-bound-repr} we prove (\ref{th:1}) and in
Section \ref{sec:c-envelope-silov} we explain how to obtain the
existence of $C^*_e(\mathcal{A})$ from Theorem \ref{th:1}.

\section{Liftings and extremals}
\label{sec:lift-extr-agler}

In Agler's approach to model theory a family is a collection of
representations of a unital algebra satisfying the first three
canonical axioms listed in the last section.  A key result in his
model theory is that an arbitrary member of a family
$\mathcal{F}_{\mathcal{A}}$ lifts to an extremal member of the family
$\mathcal{F}_{\mathcal{A}}$ (\cite{MR90a:47016}, Proposition 5.10).
We establish this result for extended families.

\begin{theorem}\label{th:2}
  If $\rho:\mathcal{A}\to \mathcal{B}(H)$ is a representation in an
  extended family $\mathcal{F}_{\mathcal{A}}$, then $\rho$ lifts to an
  extremal representation in $\mathcal{F}_{\mathcal{A}}$.
\end{theorem}

In this section we establish a preliminary version of Thorem
(\ref{th:2}) in Lemma (\ref{th:3}) below, and then indicate a proof of
the theorem based on the lemma.

Suppose the representation $\phi_\alpha:\mathcal{A}\to
\mathcal{B}(H_\alpha)$ lifts to the representation
$\phi_\beta:\mathcal{A}\to \mathcal{B}(H_\beta)$.  Then the lifting is
\textit{trivial} if $H_\alpha$ is reducing, not just invariant, for
$\phi_\beta(\mathcal{A})$.  If the only liftings of $\phi_\alpha$ are
trivial, then $\phi_\alpha$ is extremal.

If $\phi_\beta$ lifts $\phi_\alpha$, then we define a lifting
$\phi_\delta:\mathcal{A}\to \mathcal{B}(H_\delta)$ of $\phi_\beta$ to
be \textit{strongly non-trivial} with respect to $\phi_\alpha$ if
there exists an $a\in\mathcal{A}$ such that
\begin{displaymath}
  P_{H_\delta \ominus H_\beta}\phi_\delta(a)^*|_{H_\alpha}\neq 0.
\end{displaymath}
Otherwise, the lifting is \textit{weakly trivial} relative to
$\phi_\alpha$.  Finally, $\phi_\beta$ is \textit{weakly extremal}
relative to $\phi_\alpha$ if every lifting of $\phi_\beta$ is weakly
trivial relative to $\phi_\alpha$.

\begin{lemma}\label{th:3}
  Each representation $\phi_0:\mathcal{A}\to \mathcal{B}(H_0)$ in an
  extended family $\mathcal{F}_{\mathcal{A}}$ lifts to a
  representation in $\mathcal{F}_{\mathcal{A}}$ which is weakly
  extremal relative to $\phi_0$.
\end{lemma}

\begin{proof}
  The proof is by contradiction.  Accordingly, suppose $\phi_0$ does
  not lift to a weakly extremal representation relative to $\phi_0$.
  
  Let $\kappa_0$ be the cardinality of the the set of points in the
  unit sphere of ${H}_0$, $\kappa_1$ the cardinality of the set of
  elements in the unit ball of $\mathcal{A}$.  Set $\kappa =
  2^{\aleph_0 \cdot \kappa_0 \cdot \kappa_1} > \kappa_0 \cdot
  \kappa_1$.  Let $\lambda$ be the smallest ordinal greater than or
  equal to $\kappa$.  Note that there is a $C>0$ so that for
  $\pi\in\mathcal{F}_{\mathcal{A}}$, $\|\pi(a)h\| \leq C\|a\|\|h\|$
  for all $h\in H_0$ and $a\in\mathcal{A}$.
  
  Construct a chain of liftings in $\mathcal{F}_\mathcal{A}$ by
  transfinite recursion on the ordinal $\lambda$ as follows: if
  $\alpha \leq \lambda$, and $\alpha$ has a predecessor, let
  $\phi_{\alpha}$ denote a strong (with respect to ${\phi}_0$)
  nontrivial lifting of $\phi_{\alpha-1}$.  Such an lifting exists by
  the assumption that $\phi_{0}$ does not lift to a weak extremal.  If
  $\alpha$ is a limit ordinal, set $\phi_{\alpha}$ to the spanning
  representation of $\{\phi_{\delta}\}_{\delta < \alpha}$.  For any
  $h$ in the unit sphere of $H_0$ and $a$ in the unit ball of
  $\mathcal{A}$, there are at most countably many $\alpha$'s with
  predecessors where $P_{{H}_\alpha \ominus {H}_{\alpha - 1}}
  \phi_\alpha(a)^* h \neq 0$.  Since the cardinality of the set of
  ordinal numbers less than or equal to $\lambda$ and having a
  predecessor is $\kappa$, there must be an ordinal $\beta < \lambda$
  with predecessor where $P_{{H}_\beta \ominus {H}_{\beta -1}}
  \phi_\beta(a)^* h = 0$ for all $h$ in the unit sphere of $H_0$ and
  $a$ in the unit ball of $\mathcal{A}$, so that $\phi_\beta$ is a
  lifting of $\phi_{\beta-1}$ which is weakly trivial with respect to
  $\phi_0$; a contradiction, ending the proof.
\end{proof}

\begin{proof}[Proof of (\ref{th:2})]
  We use (\ref{th:3}) to prove (\ref{th:2}).  Let
  $\phi_0:\mathcal{A}\to \mathcal{B}(H_0)$ denote a given
  representation.  Lift $\phi_0$ to a representation
  $\phi_1:\mathcal{A}\to \mathcal{B}(H_1)$ which is weakly extremal
  relative to $\phi_0$.  Lift $\phi_1$ to a representation $\phi_2$
  which is weakly extremal relative to $\phi_1$.  Continuing in this
  manner, constructs a chain $\phi_j$, $j\in\mathbb{N}$, with respect
  to the partial order on liftings with the property that $\phi_j$ is
  weakly extremal relative to $\phi_{j-1}$.  The resultant spanning
  representation $\phi_\infty:\mathcal{A}\to \mathcal{B}(H_\infty)$
  lifts $\phi_0$ and it is easily checked to be extremal, since it is
  weakly extremal relative to $\phi_j$ for all $j\in\mathbb{N}$.
\end{proof}

It is not difficult to see that the restriction of an extremal to a
reducing subspace is an extremal.  Also, in (\ref{th:2}) if we were to
take the intersection of all reducing subspaces of $\phi_\infty$
containing $H_0$, we end up with the smallest reducing subspace for
$\phi_\infty$ containing $H_0$.  Restricting to this gives a minimal
extremal $\phi_e$ lifting $\phi_0$, in the sense that if $\psi$ lifts
$\phi_0$ and $\psi \leq \phi_e$, then $\psi = \phi_e$.  Of course
$\phi_e$ may still be reducible even if $\phi_0$ is irreducible.  In
addition, there may be non-isomorphic minimal extremal liftings of
$\phi_0$.
\goodbreak

\section{Dilations and Boundary Representations}
\label{sec:dilat-bound-repr}

Let $\phi_\alpha: \mathcal{A}\to \mathcal{B}(H_\alpha)$ be a
representation.  In parallel with the theory of liftings, a dilation
$\phi_\beta: \mathcal{A}\to \mathcal{B}(H_\beta)$ is termed
\textit{trivial} if $H_\alpha$ is reducing for
$\phi_\beta(\mathcal{A})$.  If the only dilations of $\phi_\alpha$ are
trivial ones, then $\phi_\alpha$ is a $\partial$-representation.

Likewise, suppose $\phi_\delta \geq \phi_\beta \geq \phi_\alpha$ in
the partial ordering for dilations, with the representations mapping
into the operators on $H_\delta$, $H_\beta$ and $H_\alpha$,
respectively.  By assumption, we can write $H_\delta = L_\delta \oplus
H_\beta \oplus M_\delta$, where $L_\delta$ and $L_\delta \oplus
H_\beta$ are invariant for $\phi_\delta$.  We say that $\phi_\alpha$
is \textit{strongly non-trivial} with respect to $\phi_\alpha$
if there exists an $a\in\mathcal{A}$ such that either
\begin{displaymath}
  P_{L_\delta} \pi_\beta(a)|_{H_\alpha} \neq 0 \qquad \text{or}\qquad
  P_{M_\delta} \pi_\beta(a)^* |_{H_\alpha} \neq 0.
\end{displaymath}
Otherwise, the dilation is said to be \textit{weakly trivial} relative
to $\phi_\alpha$.  Finally, $\phi_\beta$ is a \textit{weak
  $\partial$-representation} relative to $\phi_\alpha$ if every
lifting of $\phi_\beta$ is weakly trivial relative to $\phi_\alpha$.

\begin{lemma}\label{th:5}
  Each representation $\phi_0:\mathcal{A}\to \mathcal{B}(H_0)$ in an
  extended family $\mathcal{F}_{\mathcal{A}}$ dilates to a weak
  $\partial$-representation relative to $\phi_0$ which is also in
  $\mathcal{F}_{\mathcal{A}}$.
\end{lemma}

\begin{proof}
  The proof closely follows that of the existence of weak extremals,
  and is by contradiction.  Hence we suppose $\phi_0$ does not lift to
  a weak $\partial$-representation relative to $\phi_0$.  We define
  the ordinal $\lambda$ as in the proof of (\ref{th:3}).
  
  Construct a chain of dilations in $\mathcal{F}_\mathcal{A}$ where
  each of the representations by transfinite recursion on the ordinal
  $\lambda$ as in (\ref{th:3}): if $\alpha \leq \lambda$ and $\alpha$
  is a limit ordinal, set $\phi_{\alpha}$ to the spanning
  representation of $\{\phi_{\delta}\}_{\delta < \alpha}$ and if
  $\alpha$ has a predecessor, let $\phi_{\alpha}$ be a dilation to a
  strong (with respect to ${\phi}_0$) nontrivial dilation of
  $\phi_{\alpha-1}$, which exists by the assumption that $\phi_{0}$
  does not lift to a weak $\partial$-representation.  Then for any $h$
  in the unit sphere of $H_0$ and $a$ in the unit ball of
  $\mathcal{A}$, there are at most countably many $\alpha$'s with
  predecessors where $P_{L_\delta} \pi_\beta(a) h \neq 0$ or
  $P_{M_\delta} \pi_\beta(a)^* h \neq 0$.  The same reasoning then
  gives a representation $\phi_\beta$ in our chain dilating
  $\phi_{\beta-1}$ which is weakly trivial with respect to $\phi_0$, a
  contradiction.
\end{proof}

\begin{proof}[Proof of (\ref{th:1})]
  This now follows the proof of (\ref{th:2}).  Construct a countably
  infinite chain of representations $\{\phi_i\}$ into the bounded
  operators on Hilbert spaces $H_i$, where $\phi_i$ is a weak
  $\partial$-representation with respect to $\phi_{i-1}$ for each
  $i\in\mathbb{N}$.  Let $\phi_\infty$ denote the spanning
  representation on the Hilbert space $H_\infty$.  Since a dilation of
  a weak $\partial$-representation with respect to a representation
  $\phi$ is also a weak $\partial$-representation with respect to
  $\phi$, $\phi_\infty$ is a weak $\partial$-representation with
  respect to $\phi_i$ for all $i$.  It easily follows that
  $\phi_\infty$ is a $\partial$-representation.
\end{proof}

Minimal $\partial$-representations dilating a given representation can
be defined in the manner of minimal extremals.

\section{The $C^*$-envelope and the {\v S}ilov Ideal}
\label{sec:c-envelope-silov}

The $C^*$-envelope of the operator algebra $\mathcal{A}$, denoted
$C^*_e(\mathcal{A})$, is a $C^*$-algebra which is determined by the
property: there exists a completely isometric representation
$\gamma:\mathcal{A}\to C^*_e(\mathcal{A})$ such that $C^*(\gamma
(\mathcal{A}))=C^*_e(\mathcal{A})$ and if $\rho:\mathcal{A}\to
\mathcal{B}(H)$ is any other completely contractive representation,
then there exists an onto representation
$\pi:C^*(\rho(\mathcal{A}))\to C^*(\gamma (\mathcal{A}))$ such that
$\pi(\rho(a))=\gamma(a)$ for all $a\in\mathcal{A}$.

It is not hard to see that $C^*_e(\mathcal{A})$ is essentially unique,
for if $\rho$ also has the properties of $\gamma$, then there exists
an onto representation $\sigma:C^*(\gamma(\mathcal{A})) \to
C^*(\rho(\mathcal{A}))$ with $\sigma(\gamma(a))=\rho(a)$ for all $a\in
\mathcal{A}$.  It follows that $\sigma$ is the inverse of $\pi$ and
thus, as $C^*$-algebras, $C^*(\gamma(\mathcal{A}))$ equals
$C^*(\rho(\mathcal{A}))$.

\begin{theorem}[\cite{MR81h:46071}]\label{th:6}
  Every unital operator algebra has a $C^*$-envelope.
\end{theorem}

\begin{proof}
  A proof follows directly from (\ref{th:1}).  Viewing $\mathcal{A}$
  as a subspace of $\mathcal{B}(K)$, the inclusion mapping
  $\iota:\mathcal{A}\to \mathcal{B}(K)$ is a completely isometric
  representation and thus, according to this proposition, it dilates
  to a completely isometric representation $\gamma:\mathcal{A} \to
  \mathcal{B}(H)$ which is a $\partial$-representation.
  
  To see that $C^*(\gamma(\mathcal{A}))$ is the $C^*$-envelope,
  suppose $\psi:\mathcal{A}\to \mathcal{B}(H_\psi)$ is also completely
  isometric.  In this case $\sigma:\psi(\mathcal{A})\to
  \mathcal{B}(H)$ given by $\sigma(\psi(a))=\gamma(a)$ is completely
  contractive (and thus well-defined).  By a theorem of Arveson,
  there exists a Hilbert space $K$ containing $H$ and a representation
  $\pi:C^*(\psi(\mathcal{A}))\to \mathcal{B}(K)$ such that
  $\gamma(a)=\sigma(\psi(a))=P_H \pi(\psi(a))|_H$ (\cite{MR88h:46111},
  Cor.~6.7).  Since $a\mapsto P_H\pi(\psi(a))|_H$ is a representation
  of $\mathcal{A}$ and $\gamma$ is a $\partial$-representation, $H$
  reduces $\pi(\psi(\mathcal{A}))$.  Thus, $\sigma$ extends to an onto
  representation $C^*(\psi(\mathcal{A}))\to C^*(\gamma
  (\mathcal{A}))$.
\end{proof}

Arveson says that $\mathcal{J}$ is the {\v S}ilov boundary of the
concrete operator algebra $\mathcal{A}\subset \mathcal{B}(K)$ if
$\mathcal{J}$ contains every ideal $\mathcal{I}$ with the property
that the restriction of the quotient $q:C^*(\mathcal{A})\to
C^*(\mathcal{A})/\mathcal{I}$ to $\mathcal{A}$ is completely
isometric.  Since the inclusion of $\mathcal{A}$ into $\mathcal{B}(K)$
is completely isometric, there exists an onto representation
$\pi:C^*(\mathcal{A})\to C^*_e(\mathcal{A})=C^*(\gamma(a))$ such that
$\pi(a)=\gamma(a),$ where $\gamma$ is a representation as in
(\ref{th:6}) which generates the $C^*$-envelope of $\mathcal{A}$.  It
is left to the interested reader to verify that the kernel of $\pi$ is
the {\v S}ilov ideal of $\mathcal{A}$.

\bibliographystyle{plain}
\bibliography{b_reps}

\end{document}